\documentclass[12pt]{amsart}
\usepackage{latexsym,amsmath,amssymb,epsfig,amsthm}
\usepackage{graphicx}
\usepackage{tikz}
\usetikzlibrary{calc}
\usepackage[enableskew,vcentermath]{youngtab}
\usepackage{hyperref}
\hypersetup{colorlinks=false}
\input epsf
\newtheorem{theorem}{Theorem}[section]
\newtheorem{proposition}[theorem]{Proposition}
\newtheorem{corollary}[theorem]{Corollary}

\newtheorem{definition}[theorem]{Definition}
\newtheorem{example}[theorem]{Example}
\newtheorem{lemma}[theorem]{Lemma}

\newtheorem{question}[theorem]{Question}

\newcommand{\esd}{{\rm esd}}

\newcommand{\sd}{{\rm sd}}

\newcommand{\supp}{{\rm supp}}

\newcommand{\dD}{{\mathcal D}}
\newcommand{\eE}{{\mathcal E}}
\newcommand{\fF}{{\mathcal F}}

\newcommand{\pP}{{\mathcal P}}
\newcommand{\lL}{{\mathcal L}}
\newcommand{\qQ}{{\mathcal Q}}
\newcommand{\sS}{{\mathcal S}}
\newcommand{\uU}{{\mathcal U}}
\newcommand{\zZ}{{\mathcal Z}}

\newcommand{\RR}{{\mathbb R}}
\newcommand{\fS}{{\mathfrak S}}
\newcommand{\NN}{{\mathbb N}}
\newcommand{\ZZ}{{\mathbb Z}}

\renewcommand{\to}{\rightarrow}

\newcommand{\sm}{{\smallsetminus}}
\begin{document}
\title[Symmetric decompositions and real-rootedness]
{Symmetric decompositions, triangulations and 
real-rootedness}

\author{Christos~A.~Athanasiadis}

\address{Department of Mathematics\\
National and Kapodistrian University of Athens\\
Panepistimioupolis\\ 15784 Athens, Greece}
\email{caath@math.uoa.gr}

\author{Eleni~Tzanaki}

\address{Department of Mathematics and Applied 
Mathematics\\ University of Crete\\
70013 Heraklion, Greece}
\email{etzanaki@uoc.gr}

\date{March 4, 2021}

\begin{abstract}
Polynomials which afford nonnegative, real-rooted 
symmetric decompositions have been investigated recently 
in algebraic, enumerative and geometric combinatorics. 
Br\"and\'en and Solus have given sufficient conditions 
under which the image of a polynomial under a certain 
operator associated to barycentric subdivision has 
such a decomposition. 
This paper gives a new proof of their result which 
generalizes to subdivision operators in the setting 
of uniform triangulations of simplicial complexes, 
introduced by the first named author. Sufficient 
conditions under which these decompositions are 
also interlacing are described. Applications yield 
new classes of polynomials in geometric 
combinatorics which afford nonnegative, real-rooted 
symmetric decompositions. Some interesting questions 
in $f$-vector theory arise from this work.
\end{abstract}

\maketitle

\section{Introduction}
\label{sec:intro}

Polynomials with nonnegative coefficients and 
only real roots arise frequently in mathematics,
especially in algebra, combinatorics and 
geometry~\cite{Bra15, Sta89}. Real-rootedness 
implies strong conditions on the coefficients, 
such as unimodality and log-concavity (for missing 
definitions, see Section~\ref{sec:pre}), 
and provides a powerful method to prove these important 
properties. Real-rooted polynomials with nonnegative 
and symmetric coefficients form an especially nice 
class of polynomials. They have a property which is 
stronger than unimodality, namely $\gamma$-positivity, 
and their coefficients peak at a predictable position. 

Polynomials with nonnegative but not necessarily 
symmetric 
coefficients are amenable to techniques suitable 
for polynomials with symmetric coefficients via 
their symmetric decompositions. Every polynomial 
$p(x) \in \RR[x]$ of degree at most $n$ can be 
written uniquely in the form $p(x) = a(x) + xb(x)$
for some polynomials $a(x), b(x) \in \RR[x]$ of 
degrees at most $n$ and $n-1$, respectively, such 
that $a(x) = x^n a(1/x)$ and $b(x) = x^{n-1} 
b(1/x)$. One then hopes that these two symmetric
polynomials have nice properties and $p(x)$ is 
said to have a nonnegative, unimodal, 
$\gamma$-positive or real-rooted symmetric 
decomposition with respect to $n$ if both $a(x)$ 
and $b(x)$ have the corresponding property. 
To motivate this paper better, let us discuss two 
important examples from geometric combinatorics.

The first example comes from the theory of face 
enumeration of simplicial complexes. A convenient 
way to record the face numbers of a simplicial 
complex $\Delta$ is the $h$-polynomial, defined 
by the formula
\begin{equation}
\label{eq:hdef}
h(\Delta, x) \ = \ \sum_{i=0}^n f_{i-1} (\Delta) \, 
x^i (1-x)^{n-i} \, , 
\end{equation}
where $f_i(\Delta)$ is the number of $i$-dimensional 
faces of $\Delta$ and $n-1$ is its dimension. The 
$h$-polynomial has nonnegative coefficients if 
$\Delta$ is Cohen--Macaulay over some field and, in 
particular, if $\Delta$ triangulates a ball or 
sphere \cite[Chapter~II]{StaCCA}. If $\Delta$ 
triangulates a sphere, then $h(\Delta, x)$ has 
symmetric coefficients and its unimodality and 
$\gamma$-positivity have been major topics of 
research in the past few decades; see \cite{Ad18, APP21}
\cite[Section~3]{Ath18} \cite[Section~7.3.2]{Bra15} 
\cite{KN16, Sta89, StaCCA}. Although the 
$\gamma$-positivity of $h(\Delta, x)$ is conjectured 
to hold for all flag triangulations $\Delta$ of the 
sphere~\cite{Ga05}, no reasonable guess for when 
$h(\Delta, x)$ is real-rooted exists. It is an open 
problem to decide whether this is the case 
for barycentric subdivisions of boundary complexes of 
polytopes \cite[Question~1]{BW08}, a special class 
of flag triangulations of the sphere. An affirmative 
answer to this question has been given for simplicial 
polytopes in~\cite{BW08} and only very recently 
for cubical polytopes in~\cite{Ath20++}.

On the other hand, if $\Delta$ triangulates an 
$(n-1)$-dimensional ball, then $h(\Delta, x)$ will 
typically not have symmetric coefficients, but has 
the symmetric decomposition
\begin{equation}
\label{eq:hball-sym}
h(\Delta, x) \ = \ h(\partial \Delta, x) + 
   (h(\Delta, x) - h(\partial \Delta, x))
\end{equation}	
with respect to $n-1$, where $\partial \Delta$ 
stands for the boundary complex of $\Delta$. It 
seems natural to investigate under what conditions 
this symmetric decomposition has nice properties. 
As a consequence of \cite[Theorem~2.1]{Sta93}, 
(\ref{eq:hball-sym}) 
is nonnegative provided that no facet of $\Delta$ 
has all its vertices in $\partial \Delta$. Thus,
one expects that (\ref{eq:hball-sym}) has even 
better behavior when $\partial \Delta$ is a 
vertex-induced subcomplex of $\Delta$. Indeed, 
under this assumption, the unimodality of 
(\ref{eq:hball-sym}) follows for a large family 
of triangulations $\Delta$ of the ball from 
\cite[Theorem~46]{AY20} and one can speculate  
that (\ref{eq:hball-sym}) is also $\gamma$-positive
when $\Delta$ is flag. In fact, the latter statement 
can be shown to be equivalent to the equator 
conjecture, already posed in~\cite{CN20}. Once 
again, the real-rootedness of (\ref{eq:hball-sym})
has been much less studied. This can be deduced 
from part (b) of Theorem~\ref{thm:B-BS}, stated in 
the sequel, when $\Delta$ is the barycentric 
subdivision of a simplicial ball (see 
\cite[Section~8]{Ath20+}), but should be expected 
to hold in much more general situations. For 
instance, it seems natural to ask, in the spirit of 
\cite[Question~1]{BW08}, whether it holds for all
barycentric subdivisions of polyhedral balls.

The second example comes from Ehrhart theory.
Let $\pP \subseteq \RR^N$ be any $n$-dimensional 
convex polytope with vertices in $\ZZ^N$. The 
Ehrhart polynomial of $\pP$ \cite{BR15} \cite[Part Two]
{HiAC} \cite[Section~4.6]{StaEC1} is the unique 
polynomial $\iota(\pP; x)$ for which $\iota(\pP; m)$
is equal to the number of elements of $m \pP \cap 
\ZZ^N$ for every $m \in \NN$. The function
$h^*(\pP, x)$ defined by the formula 
\begin{equation*} %%\label{eq:h*def}
\sum_{m \ge 0} \iota(\pP; m) x^m \ = \ 
\frac{h^*(\pP, x)}{(1-x)^{n+1}}
\end{equation*}
is a well studied polynomial of degree at most $n$
with nonnegative coefficients, called the 
$h^*$-polynomial of $\pP$. Stapledon~\cite{Stap09} 
showed that $h^*(\pP, x)$ has a nonnegative symmetric 
decomposition with respect to $n$ whenever $\pP$ 
contains a lattice point in its relative interior 
(and in fact, to the best of our knowledge, the 
concept of a symmetric decomposition first appeared 
in~\cite{Stap09}). While the question of
unimodality of the $h^*$-polynomial has long been 
studied~\cite{Br16}, its $\gamma$-positivity and 
real-rootedness have been investigated more recently 
for several special classes of lattice polytopes 
\cite{BJM19, Fe20+, HJM19, Jo18, OT20, OT20+, So19} 
and the unimodality and real-rootedness of its 
symmetric decomposition have been addressed too 
\cite{BS20, Jo20+, SV13}. 

Other examples of $\gamma$-positive and real-rooted 
symmetric decompositions in enumerative combinatorics 
of geometric flavor can be found in \cite[Section~5]
{Ath18} \cite{Ath20a, Ath20b, GS20, HZ19}.

This discussion suggests that polynomials having 
nonnegative, real-rooted symmetric decompositions 
arise naturally in combinatorics and are worth of 
further study, and that developing more techniques 
to prove this property is desirable.
The starting point for this paper is the following 
theorem. Let us denote by $\RR_n[x]$ the space of 
all polynomials of degree at most $n$, with real 
coefficients.
\begin{theorem} \label{thm:B-BS} 
Let $h(x) = c_0 + c_1 x + \cdots + c_n x^n \in 
\RR_n[x]$ be a polynomial with nonnegative 
coefficients and define $\dD_n(h(x)) \in \RR_n[x]$ 
by the equation
\begin{equation} \label{eq:h-sd}
\sum_{m \ge 0} f(m) x^m \ = \ \frac{\dD_n(h(x))}
{(1-x)^{n+1}} \, ,
\end{equation}
where $f(x) = \sum_{i=0}^n c_i x^i(1+x)^{n-i}$.

\begin{itemize}
\itemsep=0pt
\item[{\rm (a)}] {\rm (\cite[Theorem~4.2]{Bra06})}
The polynomial $\dD_n(h(x))$ has only real roots.

\item[{\rm (b)}] {\rm (\cite[Theorem~2.13]{BS20})}
If the inequalities 
\begin{equation} \label{eq:c-ineq1} 
c_0 + c_1 + \cdots + c_i \ \le \ c_n + c_{n-1} + 
\cdots + c_{n-i}
\end{equation}
hold for all $0 \le i \le \lfloor n/2 \rfloor$, 
then $\dD_n(h(x))$ has a nonnegative, real-rooted 
symmetric decomposition with respect to $n$.
\end{itemize}
\end{theorem}

The map $\dD_n$ is closely related to the 
subdivision operator $\eE: \RR[x] \to \RR
[x]$; see \cite[Section~4]{Bra06} and references 
therein, \cite[Section~7.3.3]{Bra15} and 
Example~\ref{ex:sd}. It is the unique linear 
operator $\dD_n: \RR_n[x] \to \RR_n[x]$ which has 
the property that $\dD_n (h(\Delta, x)) = 
h(\sd(\Delta), x)$ for every $(n-1)$-dimensional 
simplicial complex $\Delta$, where $\sd(\Delta)$ 
is the barycentric subdivision of $\Delta$. Part
(a) of Theorem~\ref{thm:B-BS} was applied 
in~\cite{BW08} to prove that $h(\sd(\Delta), x)$ has
only real roots for every Cohen--Macaulay simplicial
complex $\Delta$. Part (b) was applied in~\cite{BS20},
among other situations, to prove that $h(\sd(\Delta), 
x)$ and $h^*(\zZ, x)$ have a nonnegative, real-rooted 
symmetric decomposition with respect to $n$ for every 
doubly Cohen--Macaulay simplicial complex $\Delta$
and every $n$-dimensional lattice zonotope $\zZ$ 
having an interior lattice point, respectively.

Since Theorem~\ref{thm:B-BS} is closely related to 
barycentric subdivision, it is natural to wonder 
whether there is a similar result which applies to 
more general types of triangulations. Indeed, part (a) 
of the theorem was generalized in~\cite{Ath20+} in 
the framework of uniform 
triangulations of simplicial complexes, of which 
barycentric subdivision is a prototypical example.
The operator $\dD_n$ is replaced there by an
operator $\dD_{\fF, n}: \RR_n[x] \to \RR_n
[x]$ which depends on a triangular array of numbers
$\fF$ and maps $h(\Delta, x)$ to the $h$-polynomial
of a triangulation of $\Delta$ for every 
$(n-1)$-dimensional simplicial complex $\Delta$, 
provided that the $f$-vector (prescribed by $\fF$) 
of the restriction of this triangulation to a face 
of $\Delta$ depends only on the dimension of that 
face (the collection of $f$-vectors of these 
restrictions is precisely the information encoded 
in $\fF$). One of the main results 
of~\cite{Ath20+} describes conditions on $\fF$ 
which guarantee that $\dD_{\fF, n}(h(x))$ has
only real roots for every polynomial $h(x) \in 
\RR_n[x]$ with nonnegative coefficients. This 
result is more general than Theorem~\ref{thm:B-BS} 
(a) and specializes to the latter and to a result 
of Jochemko~\cite{Jo18} on the real-rootedness of 
the $h^*$-polynomial in the important special cases 
of barycentric and edgewise subdivisions, 
respectively. The purpose of this paper is to prove 
that, somewhat unexpectedly,  part (b) too of 
Theorem~\ref{thm:B-BS} is valid when the operator 
$\dD_n$ is replaced by $\dD_{\fF, n}$ under the 
same assumptions on $\fF$ as those in~\cite{Ath20+}.
This provides a useful tool to address questions 
about the real-rootedness of symmetric decompositions 
in geometric combinatorics and shows that uniform 
triangulations provide a good framework to study 
this phenomenon as well. To avoid a longer discussion
in this introduction, we postpone the exact 
formulation of our main result until 
Section~\ref{sec:proof} and list some of its 
consequences instead, to demonstrate its
applicability.

The following statement is our first application. 
The operator $\uU^n_r$, defined there, 
is the operator $\dD_{\fF,n}$ associated to the $r$-fold 
edgewise subdivision. The real-rootedness of symmetric 
decompositions of polynomials of the form $\uU^n_r(h(x))$ 
was studied in the context of Ehrhart theory 
in~\cite{Jo20+}. The following theorem complements 
the results of~\cite{Jo20+}; the conclusion
of part (b) is shown in \cite[Theorem~1.1]{Jo20+} under
stronger assumptions on $h(x)$ (but for a possibly 
larger range of values of $r$).
\begin{theorem} \label{thm:main-esdr} 
Let $h(x) = c_0 + c_1 x + \cdots + c_n x^n \in \RR_n[x]$ 
be a polynomial with nonnegative coefficients. Given a 
positive integer $r$, define $\uU^n_r(h(x)) \in \RR_n[x]$ 
by the formula
\begin{equation*} %%\label{eq:def-Unr-intro}
\frac{\uU^n_r(h(x))}{(1-x)^n} \ = \ \sum_{m \ge 0} 
a_{rm} x^m \ \ \ \ \ \text{if} \ \ \ \ 
\frac{h(x)}{(1-x)^n} \ = \ \sum_{m \ge 0} a_m x^m.
\end{equation*}

\begin{itemize}
\itemsep=0pt
\item[{\rm (a)}]
If the inequalities (\ref{eq:c-ineq1}) 
hold for $0 \le i \le \lfloor n/2 \rfloor$, then 
$\uU^n_r(h(x))$ has a nonnegative, real-rooted 
symmetric decomposition with respect to $n$ for 
every $r \ge n$.

\item[{\rm (b)}] 
If $c_n = 0$ and the inequalities 
\begin{equation} \label{eq:c-ineq2} 
c_0 + c_1 + \cdots + c_i \ \ge \ c_{n-1} + c_{n-2} 
+ \cdots + c_{n-i}
\end{equation}
hold for all $1 \le i \le \lfloor n/2 \rfloor$,
then $\uU^n_r(h(x))$ has a nonnegative, real-rooted 
symmetric decomposition with respect to $n-1$ for 
every $r \ge n$.
\end{itemize}
\end{theorem}

Our second application generalizes 
Theorem~\ref{thm:B-BS}. Indeed, the operator 
$\dD_{n,r}$, defined in the 
following statement, reduces to $\dD_n$ for $r=1$; 
it coincides with the operator $\dD_{\fF,n}$ defined 
by a generalization of barycentric subdivision, termed 
as the $r$-colored barycentric subdivision 
in~\cite{Ath20+}. Part (a) coincides with 
\cite[Proposition~7.5]{Ath20+}.

\begin{theorem} \label{thm:main-sdr} 
Let $h(x) = c_0 + c_1 x + \cdots + c_n x^n \in 
\RR_n[x]$ be a polynomial with nonnegative 
coefficients. Given a positive integer $r$, define 
$\dD_{n,r}(h(x)) \in \RR_n[x]$ by the equation
\begin{equation*} %%\label{eq:h-sdr}
f(0) + \sum_{m \ge 1} \left( f(rm) - f(rm-1) 
\right) x^m \ = \ \frac{\dD_{n,r}(h(x))} 
{(1-x)^n} \, ,
\end{equation*}
where $f(x) = \sum_{i=0}^n c_i x^i (1+x)^{n-i}$.

\begin{itemize}
\itemsep=0pt
\item[{\rm (a)}]
The polynomial $\dD_{n,r}(h(x))$ has only real 
roots.

\item[{\rm (b)}] 
If the inequalities (\ref{eq:c-ineq1}) hold for 
all $0 \le i \le \lfloor n/2 \rfloor$, then 
$\dD_{n,r}(h(x))$ has a nonnegative, real-rooted 
symmetric decomposition with respect to $n$.
\end{itemize}
\end{theorem}

Applications of Theorems~\ref{thm:main-esdr} 
and~\ref{thm:main-sdr} to the $h$-polynomials
of $r$-fold edgewise subdivisions and $r$-colored 
barycentric subdivisions of doubly Cohen--Macaulay 
simplicial complexes and triangulations of balls are 
given in Section~\ref{sec:apps}. The symmetric 
decomposition (\ref{eq:hball-sym}), in particular, 
is shown there to be real-rooted for new classes of 
triangulations of the ball.

A lattice zonotope is 
defined as the Minkowski sum of finitely many line 
segments in $\RR^N$, whose vertices lie in $\ZZ^N$. 
The following statement can be deduced from 
Theorem~\ref{thm:main-sdr} (see Section~\ref{sec:apps}).
Recall that $\iota(\zZ; x)$ stands for the Ehrhart 
polynomial of a lattice polytope $\zZ$ and note that 
the polynomial $h_r^*(\zZ,x)$, defined 
by (\ref{eq:r-zono}), reduces to the $h^*$-polynomial 
$h^*(\zZ,x)$ for $r=1$. Thus, parts (a) and (b) of the
following statement generalize the main results of 
\cite{BJM19} and \cite[Section~4]{BS20}, respectively. 
In the notation of Theorem~\ref{thm:main-esdr}, we 
have $h_r^*(\zZ,x) = \uU^n_r(h^*(\zZ,x))$ for every 
$r \ge 1$. 
\begin{corollary} \label{cor:zono}
Let $\zZ$ be an $n$-dimensional lattice zonotope and 
for positive integers $r$, define $h_r^*(\zZ,x) \in
\RR_n[x]$ by the equation 
\begin{equation} \label{eq:r-zono}
1 + \sum_{m \ge 1} 
\left( \iota(\zZ; rm) - \iota(\zZ; rm-1) \right) x^m 
\ = \ \frac{h_r^*(\zZ,x)}{(1-x)^n} \, .
\end{equation}

\begin{itemize}
\itemsep=0pt
\item[{\rm (a)}]
The polynomial $h_r^*(\zZ,x)$ has only real roots
for every $r \ge 1$.

\item[{\rm (b)}] 
If $\zZ$ has a lattice point in its relative 
interior, then $h_r^*(\zZ,x)$ has a nonnegative, 
real-rooted symmetric decomposition with respect to 
$n$ for every $r \ge 1$.
\end{itemize}
\end{corollary}

Our final application identifies a class of doubly 
Cohen--Macaulay simplicial complexes, namely that 
of one-coskeleta of Cohen--Macaulay simplicial 
complexes, whose uniform triangulations have 
$h$-polynomials with especially nice symmetric 
decompositions. This result, stated here for
edgewise and $r$-colored barycentric subdivisions, 
is new even for barycentric subdivisions. We recall 
that a nonnegative, real-rooted symmetric 
decomposition $p(x) = a(x) + xb(x)$ is said to be 
interlacing if $a(x)$ is interlaced by $b(x)$ 
(see \cite[Theorem~2.6]{BS20} for a number of 
equivalent conditions).
\begin{theorem} \label{thm:skel-intro} 
Let $\Gamma$ be any $n$-dimensional simplicial 
complex with nonnegative $h$-vector and let $\Delta$ 
be the $(n-1)$-dimensional skeleton of $\Gamma$. 

\begin{itemize}
\itemsep=0pt
\item[{\rm (a)}]
The polynomial $\uU^n_r(h(\Delta, x))$ has a 
nonnegative, real-rooted and interlacing symmetric 
decomposition with respect to $n$ for every $r \ge n$.
Moreover, $\uU^n_r(h(\Delta, x))$ interlaces 
$\uU^{n+1}_r(h(\Gamma, x))$ for every $r \ge n+1$.

\item[{\rm (b)}] 
The polynomial $\dD_{n,r}(h(\Delta, x))$ has a 
nonnegative, real-rooted and interlacing symmetric 
decomposition with respect to $n$ for every $r \ge 1$.
Moreover, $\dD_{n,r}(h(\Delta, x))$ interlaces 
$\dD_{n+1,r}(h(\Gamma, x))$.
\end{itemize}
\end{theorem}

The proof of Theorem~\ref{thm:B-BS}, given 
in~\cite{BS20}, uses a lot of technical properties 
of the subdivision operators $\dD_n$ 
and $\eE$. The proof of our more general theorem, 
given in Section~\ref{sec:proof}, involves no such 
technicalities and essentially uses only the universal 
recurrence for the polynomials $\dD_{\fF, n}(x^k)$
\cite[Lemma~6.3]{Ath20+} and basic facts about 
real-rooted polynomials.

The structure of the remainder of this 
paper is as follows. Section~\ref{sec:pre} fixes 
notation and recalls useful definitions and 
facts about simplicial complexes and real-rooted 
polynomials. Section~\ref{sec:uniform} reviews the 
basics of the enumerative theory of uniform 
triangulations~\cite{Ath20+}. The main result 
(Theorem~\ref{thm:main}) of this paper is stated and 
proven in Section~\ref{sec:proof}. The sufficient 
conditions provided for the real-rooted symmetric 
decompositions, considered there, to be interlacing 
lead to new inequalities that the $h$-vector of a 
Cohen--Macaulay simplicial complex may or may not 
satisfy (see Corollary~\ref{cor:main}) and raise 
questions in $f$-vector theory which are of 
independent interest (see Section~\ref{sec:rem}). 
Theorems~\ref{thm:main-esdr} and~\ref{thm:main-sdr} 
are deduced from Theorem~\ref{thm:main} in 
Section~\ref{sec:apps} and some of their own 
consequences (see Corollaries~\ref{cor:esdr} 
and~~\ref{cor:sdr}) are discussed there. For the 
proof of Theorem~\ref{thm:main-sdr}, one needs to 
verify that the operator $\dD_{n,r}$ satisfies the 
crucial conditions of \cite[Theorem~6.1]{Ath20+}, 
a problem that was left open in 
\cite[Section~7]{Ath20+}. This nontrivial fact 
requires special treatment and is proven in 
Section~\ref{sec:r-colored}. 
Theorem~\ref{thm:skel-intro} is stated more 
generally, in the setting of uniform triangulations, 
and proven in Section~\ref{sec:skeleta}. 
Section~\ref{sec:rem} concludes with remarks and
questions that are raised by this work.

\section{Preliminaries}
\label{sec:pre}

This section fixes notation and explains background 
and terminology on real polynomials, simplicial 
complexes and their triangulations which will be 
useful in the sequel.

\subsection{Polynomials}
We recall that $\RR_n[x]$ stands for the space of 
polynomials of degree at most $n$ with real coefficients.
A polynomial $p(x) = a_0 + a_1 x + \cdots + a_n x^n 
\in \RR_n[x]$ is called
\begin{itemize}
\item[$\bullet$] 
  \emph{symmetric}, with center of symmetry $n/2$, if 
	$a_i = a_{n-i}$ for all $0 \le i \le n$,
\item[$\bullet$] 
  \emph{unimodal}, with a peak at position $k$, if $a_0 
	\le a_1 \le \cdots \le a_k \ge a_{k+1} \ge \cdots \ge 
	a_n$,
\item[$\bullet$] 
  \emph{$\gamma$-positive}, with center of symmetry $n/2$,
	if $p(x) = \sum_{i=0}^{\lfloor n/2 \rfloor} \gamma_i 
	x^i (1+x)^{n-2i}$ for some nonnegative reals $\gamma_0, 
  \gamma_1,\dots,\gamma_{\lfloor n/2 \rfloor}$,
\item[$\bullet$] 
   \emph{log-concave}, if $a^2_i \ge a_{i-1}a_{i+1}$ for 
	 all $1 \le i < n$,
\item[$\bullet$] 
  \emph{real-rooted}, if every root of $p(x)$ is real, or 
	$p(x) = 0$.
\end{itemize}
Every $\gamma$-positive polynomial is symmetric and 
unimodal and every real-rooted and symmetric polynomial
with nonnegative coefficients is $\gamma$-positive; see
\cite{Ath18, Bra15, Sta89} for more information on the 
connections among these concepts.  

A real-rooted polynomial $p(x)$, with 
roots $\alpha_1 \ge \alpha_2 \ge \cdots$, is said to 
\emph{interlace} a real-rooted polynomial $q(x)$, with 
roots $\beta_1 \ge \beta_2 \ge \cdots$, if
\[ \cdots \le \alpha_2 \le \beta_2 \le \alpha_1 \le
   \beta_1. \]
By convention, the zero polynomial interlaces and is 
interlaced by every real-rooted polynomial. 
A sequence $(p_0(x), p_1(x),\dots,p_m(x))$ of 
real-rooted polynomials with nonnegative coefficients 
is called \emph{interlacing} if $p_i(x)$ interlaces 
$p_j(x)$ for $0 \le i < j \le m$. The importance of
this concept for us comes from the fact that every 
nonnegative linear combination $p(x)$ of $p_0(x), 
p_1(x),\dots,p_m(x)$ is then real-rooted; moroever, 
$p(x)$ interlaces $p_m(x)$ and is interlaced by 
$p_0(x)$. 

A standard way to produce interlacing sequences in 
combinatorics is the following. Suppose that $p_0(x), 
p_1(x),\dots,p_m(x)$ are real-rooted polynomials 
with nonnegative coefficients and set 
\begin{equation} \label{eq:recipe}
  q_k(x) \ = \ x \sum_{i=0}^{k-1} p_i(x) \, + \, 
  \sum_{i=k}^m p_i(x) 
\end{equation}
for $k \in \{0, 1,\dots,m+1\}$. Then, if the sequence 
$(p_0(x), p_1(x),\dots,p_m(x))$ is interlacing, so is 
$(q_0(x), q_1(x),\dots,q_{m+1}(x))$; see 
\cite[Corollary~8.7]{Bra15} for a more general 
statement. For an extensive treatment of real-rooted 
polynomials and the theory of interlacing, 
see~\cite{Fi06}.

The $k$th Veronese $r$-section operator is defined 
on polynomials, or formal power series, by the 
formula
\[ \sS^r_k \left( \, \sum_{n \ge 0} a_nx^n \right) 
   \, = \ \sum_{n \ge 0} a_{rn+k} x^n . \]
We note that
\begin{equation} \label{eq:Sx}
\sS^r_j (x^i f(x)) \ = \ \begin{cases}
  \sS^r_{j-i} (f(x)), & \text{if $i \le j$} \\
  x \sS^r_{r-i+j} (f(x)), & \text{if $i > j$} 
\end{cases} \end{equation}
for $i, j \in \{0, 1,\dots,r-1\}$.

\subsection{Simplicial complexes}
We assume familiarity with basic notions from  
algebraic, enumerative and topological combinatorics 
on simplicial complexes; excellent resources on these 
topics are \cite{Bj95, HiAC, StaCCA}. All simplicial 
complexes considered here will be abstract and finite. 
Following~\cite{Ath20+}, we denote by $\sigma_n$ the 
abstract simplex $2^V$ on an $n$-element vertex set
$V$. 

For the remainder of this section, $\Delta$ will
be an $(n-1)$-dimensional simplicial complex. The 
sequence $h(\Delta) := (h_0(\Delta), 
h_1(\Delta),\dots,h_n(\Delta))$ of coefficients of 
the $h$-polynomial $h(\Delta, x) = \sum_{i=0}^n h_i
(\Delta)x^i$, already defined in the introduction 
by Equation~(\ref{eq:hdef}), is called the 
\emph{$h$-vector} of $\Delta$. As mentioned there,  
$h(\Delta)$ has nonnegative entries whenever 
$\Delta$ is Cohen--Macaulay (over some field). 
We note that $h_0(\Delta) = 1$ and $h_n
(\Delta) = (-1)^{n-1} \widetilde{\chi}(\Delta)$, 
where $\widetilde{\chi}(\Delta)$ is the reduced 
Euler characteristic of $\Delta$; in particular,
$h_n(\Delta) = 0$ if the geometric realization of
$\Delta$ is contractible. We will be interested in
simplicial complexes which satisfy the inequalities
\begin{equation} \label{eq:ineqCM*} 
h_0(\Delta) + h_1(\Delta) + \cdots + h_i(\Delta) 
\ \le \ h_n(\Delta) + h_{n-1}(\Delta) + \cdots + 
h_{n-i}(\Delta) 
\end{equation}
for $0 \le i \le \lfloor n/2 \rfloor$ 
(equivalently, for $0 \le i \le n$). Doubly 
Cohen--Macaulay simplicial complexes 
\cite[Section~III.3]{StaCCA} have this property. A 
larger family of simplicial complexes which satisfy
inequalities~(\ref{eq:ineqCM*}) was introduced and 
studied in~\cite{MM16} under the name 
\emph{uniformly Cohen--Macaulay} simplicial 
complexes. Here we will use the term 
\emph{Cohen--Macaulay*} simplicial complex instead,
to avoid confusion with our terminology ``uniform 
triangulation". Thus, $\Delta$ is Cohen--Macaulay*
if $\Delta$ and the simplicial complexes obtained 
from it by removing any (single) facet of $\Delta$ are 
Cohen--Macaulay of dimension $n-1$. Every doubly 
Cohen--Macaulay simplicial complex is Cohen--Macaulay*
(see \cite[Proposition~2.8]{MM16}) and every 
Cohen--Macaulay* simplicial complex 
satisfies~(\ref{eq:ineqCM*}) for
all $i$ \cite[Proposition~2.7]{MM16}.

By the term \emph{triangulation} of $\Delta$ we 
will always mean a geometric triangulation. Thus, 
a simplicial complex $\Delta'$ is a triangulation 
of $\Delta$ if there exists a geometric realization 
of $\Delta'$ which geometrically subdivides one for 
$\Delta$. 

Barycentric and edgewise subdivisions are 
important triangulations of $\Delta$. The 
\emph{barycentric subdivision} of $\Delta$, denoted 
by $\sd(\Delta)$, is defined as the simplicial complex 
of all chains of nonempty faces of $\Delta$. The 
edgewise subdivision depends on a positive integer 
$r$ and a linear ordering of the vertex set 
$V(\Delta)$ of $\Delta$ (although its face numbers 
are independent of the latter). Given such an ordering 
$v_1, v_2,\dots,v_m$, we denote by $V_r(\Delta)$ the set 
of maps $f: V(\Delta) \to \NN$ such that $\supp(f) \in 
\Delta$ and $f(v_1) + f(v_2) + \cdots + f(v_m) = r$, 
where $\supp(f)$ is the set of all $v \in V(\Delta)$ 
for which $f(v) \ne 0$. For $f \in V_r(\Delta)$, we let 
$\iota(f): V(\Delta) \to \NN$ be the map defined by 
setting $\iota(f)(v_j) = f(v_1) + f(v_2) + \cdots + 
f(v_j)$ for $j \in \{1, 2,\dots,m\}$. The 
\emph{$r$-fold edgewise subdivision} of $\Delta$, 
denoted by $\esd_r(\Delta)$, is the simplicial 
complex on the vertex set $V_r(\Delta)$ of which a 
set $E \subseteq V_r(\Delta)$ is a face if the 
following two conditions are satisfied:

\begin{itemize}
\itemsep=0pt
\item[$\bullet$]
$\bigcup_{f \in E} \, \supp(f) \in \Delta$ and

\item[$\bullet$]
$\iota(f) - \iota(g) \in \{0, 1\}^{V(\Delta)}$, or 
$\iota(g) - \iota(f) \in \{0, 1\}^{V(\Delta)}$, for 
all $f, g \in E$.
\end{itemize}

The simplicial complexes $\sd(\Delta)$ and 
$\esd_r(\Delta)$ can be realized as triangulations 
of $\Delta$. This is well known for the former; for 
the latter, see \cite[Section~5]{Ath20b} and 
references therein. 

A simplicial complex $\Delta$ is called \emph{flag}  
if every clique in the one-skeleton of $\Delta$ is 
a face of $\Delta$; see \cite[Section~3]{Ath18}
\cite[Section~5.2]{KN16} \cite[Section~III.4]{StaCCA} 
for information about this very interesting class 
of simplicial complexes.

\section{Uniform triangulations and subdivision 
operators} \label{sec:uniform}

This section summarizes the background on uniform
triangulations of simplicial complexes and their 
associated subdivision operators~\cite{Ath20+}
which are necessary in order to state and prove our 
main results.

An \emph{$f$-triangle} of size 
$d \in \NN \cup \{ \infty \}$ is simply a triangular 
array $\fF = (f_\fF(i,j))_{0 \le i \le j \le d}$ of 
nonnegative integers (where $i, j$ are finite numbers). 
A triangulation $\Delta'$ of a simplicial complex 
$\Delta$ of dimension less than $d$ is called 
\emph{$\fF$-uniform} if for all $0 \le i \le j \le d$, 
the restriction of $\Delta'$ to any face of $\Delta$ 
of dimension $j-1$ has exactly $f_\fF(i,j)$ faces of 
dimension $i-1$. We say that $\fF$ is \emph{feasible} 
if every simplex of dimension less than $d$ has an 
$\fF$-uniform triangulation. The barycentric 
subdivision $\sd(\Delta)$ and the $r$-fold edgewise 
subdivision $\esd_r(\Delta)$ are prototypical examples
of uniform triangulations of $\Delta$.

For every $f$-triangle 
$\fF$ of size $d$, there exist linear operators 

\medskip
\begin{center}
\begin{tabular}{rl}
$\eE_\fF$ : & $\RR_d[x] \to \RR_d[x]$ \\
$\dD_{\fF,n}$ : & $\RR_n[x] \to \RR_n[x], \ 
\text{{\rm for}} \ n \in \{0, 1,\dots,d\} \sm 
                 \{ \infty \}$
\end{tabular}
\end{center}

\medskip
\noindent
such that $f(\Delta',x) = \eE_\fF(f(\Delta,x))$ and 
$h(\Delta',x) = \dD_{\fF,n}(h(\Delta,x))$ for every 
simplicial complex $\Delta$ of dimension $n-1$, every
$\fF$-uniform triangulation $\Delta'$ of $\Delta$ and 
all finite $n \le d$. Thus, setting $p_{\fF,n,k}(x) := 
\dD_{\fF,n}(x^k)$ for $k \in \{0, 1,\dots,n\}$, we 
have 
\begin{equation} \label{eq:DF}
\dD_{\fF,n}(h(x)) \ = \ \sum_{k=0}^n c_k 
                        p_{\fF,n,k}(x) 
\end{equation}
for every polynomial $h(x) = c_0 + c_1 x + \cdots + 
c_n x^n \in \RR_n[x]$ and 
\begin{equation} \label{eq:hF}
h(\Delta',x) \ = \ \sum_{k=0}^n h_k(\Delta)  
                        p_{\fF,n,k}(x) 
\end{equation}
for every $\fF$-uniform triangulation $\Delta'$ of 
an $(n-1)$-dimensional simplicial complex $\Delta$.
Following the notation of~\cite{Ath20+}, we
write $h_\fF(\Delta, x)$ for the right-hand side of
Equation~(\ref{eq:hF}), so that $h(\Delta', x) = 
h_\fF(\Delta,x)$ for every $\fF$-uniform triangulation 
$\Delta'$ of $\Delta$.
\begin{example} \label{ex:sd} \rm 
By \cite[Section~5]{Ath20+} we have $\eE_\fF(x^n) = 
\sum_{k=0}^n f_\fF^\circ(k,n) x^k$ for every $n \le 
d$, where $f_\fF^\circ(k,n)$ is the number of interior 
$(k-1)$-dimensional faces of any $\fF$-uniform 
triangulation of the simplex $\sigma_n$.

In particular, for the $f$-triangle associated to 
barycentric subdivision we have $\eE_\fF(x^n) = 
\sum_{k=0}^n k! S(n,k) x^k$ for every $n \in \NN$, 
where $S(n,k)$ are the Stirling numbers of the second 
kind. Thus, $\eE_\fF: \RR[x] \to \RR[x]$ coincides 
with the subdivision operator $\eE: \RR[x] \to \RR[x]$ 
of \cite[Section~7.3.3]{Bra15}, mentioned in the 
introduction, and $\dD_{\fF,n}: \RR_n[x] \to \RR_n[x]$ 
coincides (see, for instance, \cite[Lemma~2.7]{BS20}) 
with the operator $\dD_n: \RR_n[x] \to \RR_n[x]$ of 
Theorem~\ref{thm:B-BS}.
\qed
\end{example} 

The polynomial $p_{\fF,n,k}(x)$ was shown to have 
nonnegative coefficients \cite[Theorem~4.1]{Ath20+} 
for every $k \in \{0, 1,\dots,n\}$ and every feasible 
$f$-triangle $\fF$ of size at least $n$. Following the
notation of \cite{Ath20+} \cite[Section~5]{ABK20+}, 
we also set 
\[ p_{\fF,n-1,n} (x) \ = \ \theta_\fF(\sigma_n,x) \
   := \ h_\fF(\sigma_n, x) - h_\fF(\partial \sigma_n, 
	 x) \]
and we consider the sequences
\begin{eqnarray*}
\pP_{\fF,n} & := & (p_{\fF,n-1,0}(x), 
               p_{\fF,n-1,1}(x),\dots,p_{\fF,n-1,n}(x)) \\ 
\qQ_{\fF,n} & := & (p_{\fF,n,0}(x), 
                    p_{\fF,n,1}(x),\dots,p_{\fF,n,n}(x)).
\end{eqnarray*}

The polynomial $\theta_\fF(\sigma_n,x)$ does not always 
have nonnegative coefficients. This is the case under 
some mild assumptions on the triangulation which defines
$\fF$; see \cite[Remark~6.1~(b)]{Ath20+} or our 
discussion in the introduction. Let us introduce the 
following useful terminology. 
\begin{definition} \label{def:lace-property}
We say that a feasible $f$-triangle $\fF$ of size at 
least $n$ has the interlacing property with respect to 
$n$ if $\qQ_{\fF,m}$ is an interlacing sequence for every 
$m \in \{0,1,\dots,n\}$, and that $\fF$ has the strong 
interlacing property with respect to $n$, if the 
following conditions hold: 

\begin{itemize}
\itemsep=0pt
\item[{\rm (i)}]
$h_\fF(\sigma_m, x)$ is a real-rooted polynomial for 
all $2 \le m < n$.

\item[{\rm (ii)}]
$\theta_\fF(\sigma_m, x)$ is either identically zero, 
or a real-rooted polynomial of degree $m-1$ with 
nonnegative coefficients which is interlaced by $h_\fF
(\sigma_{m-1}, x)$, for all $2 \le m \le n$.
\end{itemize}
\end{definition} 

We also say that a feasible $f$-triangle of infinite 
size has the (strong) interlacing property, if it does 
so with respect to every $n \in \NN$. 

The proof of \cite[Theorem~6.1]{Ath20+} shows that if 
$\fF$ has the strong interlacing property with respect 
to $n$, then $\qQ_{\fF,n}$ and $\pP_{\fF,m}$ for
$m \le n$ are interlacing sequences (that hopefully 
explains our terminology). Thus, given also that 
$p_{\fF,n,0}(x) = h_\fF(\sigma_n, x)$ and $p_{\fF,n,n}
(x) = x^n h_\fF(\sigma_n, 1/x)$ (see our 
discussion in the sequel), the following statement is 
included in the results of~\cite{Ath20+}.
\begin{theorem} \label{thm:Ath20+} 
{\rm (\cite{Ath20+})}
Let $\fF$ be any feasible $f$-triangle of size at least 
$n$ which has the strong interlacing property with 
respect to $n$. Then, $\fF$ has the interlacing property 
with respect to $n$. In particular: 

\begin{itemize}
\itemsep=0pt
\item[$\bullet$]
$\dD_{\fF,n}(h(x))$ is real-rooted, is
interlaced by $h_\fF(\sigma_n, x)$ and it interlaces 
$x^n h_\fF(\sigma_n, 1/x)$ for every polynomial 
$h(x) \in \RR_n[x]$ with nonnegative coefficients.

\item[$\bullet$]
$h_\fF(\Delta, x)$ is real-rooted, is interlaced by 
$h_\fF(\sigma_n, x)$ and it interlaces 
$x^n h_\fF(\sigma_n, 1/x)$ for every $(n-1)$-dimensional 
simplicial complex $\Delta$ with nonnegative 
$h$-vector.
\end{itemize}
\end{theorem}

The crucial strong interlacing property is especially
easy to verify for the barycentric subdivision 
\cite[Example~7.1]{Ath20+}, since then $\theta_\fF
(\sigma_n, x) = 0$ for every $n \in \NN$. It was also 
verified for the $r$-fold edgewise subdivision when 
$r \ge n$ and for certain triangulations interpolating
between barycentric and edgewise subdivisions
\cite[Section~7]{Ath20+}. Moreover, it was conjectured 
to hold for the antiprism triangulation 
\cite[Section~5]{ABK20+}, in which case only the claim 
about interlacing in condition (ii) of 
Definition~\ref{def:lace-property} is open. We will
also verify the strong interlacing property for the 
$r$-colored barycentric subdivison for every positive 
integer $r$ in Section~\ref{sec:r-colored} and will 
deduce from that and Theorem~\ref{thm:main} many of 
the results stated in the introduction.

The following proposition collects some useful 
properties of the polynomials $p_{\fF,n,k} (x)$. 
\begin{proposition} \label{prop:Ath20+}
{\rm (\cite{Ath20+})}
For every feasible $f$-triangle $\fF$ of size at 
least $n$ we have: 
\begin{eqnarray}
x^n p_{\fF,n,k} (1/x) & = & p_{\fF,n,n-k} (x), 
\ \text{{\rm for}} \ k \in \{0, 1,\dots,n\}
\label{eq:pnk-sym} \\ & & \nonumber \\
x^n p_{\fF,n-1,n}(1/x) & = & p_{\fF,n-1,n}(x),
\label{eq:theta-sym} \\ & & \nonumber \\
p_{\fF,n,k}(x) & = & x \sum_{i=0}^{k-1} 
p_{\fF,n-1,i}(x) \, + \, \sum_{i=k}^n p_{\fF,n-1,i}(x), 
\ \text{{\rm for}} \ k \in \{0, 1,\dots,n\}. 
\label{eq:pnk-rec-long}
\end{eqnarray}
\end{proposition}

As special cases of Equation~(\ref{eq:hF}), and 
applying the recurrence~(\ref{eq:pnk-rec-long}) for
$k=0$, we also have
\begin{eqnarray} 
h_\fF(\sigma_n,x) & = & p_{\fF,n,0} (x) \ = \ 
\sum_{k=0}^n p_{\fF,n-1,k}(x) , \nonumber \\ %%\label{eq:pn0}
h_\fF(\partial \sigma_n, x) & = & \sum_{k=0}^{n-1} 
p_{\fF,n-1,k}(x) . 
\label{eq:pnk-sum}
\end{eqnarray}

The following lemma will be useful in 
Section~\ref{sec:proof}.
\begin{lemma} \label{lem:DF-recip}
Let $h(x) \in \RR_{n-1}[x]$. 

\begin{itemize}
\itemsep=0pt
\item[{\rm (a)}]
$\dD_{\fF,n}(x^n h(1/x)) = x^n 
     \dD_{\fF,n}(h(x))_{x \, \mapsto 1/x}$.

\item[{\rm (b)}] 
The symmetric decomposition of $\dD_{\fF,n}(h(x))$ 
with respect to $n-1$ is nonnegative and real-rooted 
(respectively, nonnegative, real-rooted and 
interlacing), if and only if so is the symmetric 
decomposition of $\dD_{\fF,n}(x^n h(1/x))$ with 
respect to $n$. 
\end{itemize}
\end{lemma}

\begin{proof}
Part (a) follows from Equation~(\ref{eq:DF}) and 
the symmetry property~(\ref{eq:pnk-sym}). It also 
implies that for the symmetric decompositions 
$\dD_{\fF,n}(h(x)) = a(x) + x b(x)$ and 
$\dD_{\fF,n}(x^n h(1/x)) = \tilde{a}(x) + x 
\tilde{b}(x)$ of $\dD_{\fF,n}(h(x))$ and $\dD_{\fF,n}
(x^n h(1/x))$ with respect to $n-1$ and $n$, 
respectively, one has $\tilde{a}(x) = x^{n-1} 
b(1/x)$ and $\tilde{b}(x) = x^{n-1} a(1/x)$. Part 
(b) follows from these facts.
\end{proof}

\section{The main theorem}
\label{sec:proof}

This section states and proves the main results of 
this paper, using only the theory of 
Section~\ref{sec:uniform} and basic facts about 
real-rooted polynomials. Throughout it, $\fF$ stands 
for a feasible $f$-triangle of size at least $n$. 
\begin{theorem} \label{thm:main}
Let $\fF$ be a feasible $f$-triangle which has the 
strong interlacing property with respect to $n$. 
Let $h(x) = c_0 + c_1 x + \cdots + c_n x^n \in \RR_n
[x]$ be a polynomial with nonnegative coefficients.

\begin{itemize}
\itemsep=0pt
\item[{\rm (a)}] 
If the inequalities~(\ref{eq:c-ineq1}) hold for 
$0 \le i \le \lfloor n/2 \rfloor$, then $\dD_{\fF,n}
(h(x))$ has a nonnegative, real-rooted symmetric 
decomposition with respect to $n$.

If, additionally, $c_i c_{n-i-1} \le c_{i+1} c_{n-i}$ 
for all $0 \le i \le n-1$, then this decomposition 
is also interlacing.

\item[{\rm (b)}] 
If $c_n = 0$ and the inequalities~(\ref{eq:c-ineq2}) 
hold for $1 \le i \le \lfloor n/2 \rfloor$, then 
$\dD_{\fF,n}(h(x))$ has a nonnegative, real-rooted 
symmetric decomposition with respect to $n-1$.

If, additionally, $c_i c_{n-i-1} \ge c_{i+1} c_{n-i}$ 
for all $1 \le i < n-1$, then this decomposition is 
also interlacing.
\end{itemize}
\end{theorem}

The proof is based on the following lemma.
\begin{lemma} \label{lem:sym-dec}
For every $h(x) = c_0 + c_1 x + \cdots + c_n x^n \in 
\RR_n[x]$, the symmetric decomposition $\dD_{\fF,n}
(h(x)) = a(x) + xb(x)$ of $\dD_{\fF,n}(h(x))$ with 
respect to $n$ is given by
\begin{eqnarray} \label{eq:a} 
a(x) & = & \left( c_0 + c_1 + \cdots + c_n \right) 
           p_{\fF,n-1,n}(x) \, + \\ \nonumber
& & \sum_{i=0}^{n-1} \left(
c_0 + c_1 + \cdots + c_i + 
		(c_0 + c_1 + \cdots + c_{n-i-1}) x \right) 
		 p_{\fF,n-1,i}(x) \\ \nonumber & & \\ \label{eq:b} 
b(x) & = & \sum_{i=0}^{n-1} \left(c_n + c_{n-1} + 
\cdots + c_{n-i} - c_0 - c_1 - \cdots - c_i \right) 
p_{\fF,n-1,i}(x) . 
\end{eqnarray}
\end{lemma}

\begin{proof}
Let $a(x)$ and $b(x)$ be defined by~(\ref{eq:a}) 
and~(\ref{eq:b}), respectively. Properties 
(\ref{eq:pnk-sym}) and (\ref{eq:theta-sym}) of the 
$p_{\fF,n-1,k}(x)$ directly imply that $x^n a(1/x) = 
a(x)$ and $x^{n-1} b(1/x) = b(x)$. Moreover, using 
recurrence (\ref{eq:pnk-rec-long}), we get
\[ \dD_{\fF,n}(h(x)) \ = \
\sum_{k=0}^n c_k p_{\fF,n,k}(x) \ = \ 
\sum_{k=0}^n c_k \left( \, x \sum_{i=0}^{k-1} 
p_{\fF,n-1,i}(x) \, + \, \sum_{i=k}^n p_{\fF,n-1,i}
(x) \right) . \]
Changing the order of summation gives 
\begin{eqnarray*} 
\dD_{\fF,n}(h(x)) & = & 
\sum_{i=0}^n \left( \, \sum_{k > i} c_k x \, + \, 
\sum_{k=0}^i c_k \right) p_{\fF,n-1,i}(x) \ = \ 
\left( c_0 + c_1 + \cdots + c_n \right) 
p_{\fF,n-1,n}(x) \, + \\ & &
\sum_{i=0}^{n-1} \left( c_0 + c_1 + \cdots + c_i + 
c_{i+1}x + \cdots + c_n x \right) p_{\fF,n-1,i}(x) 
\ = \ a(x) + xb(x)  
\end{eqnarray*}

\medskip
\noindent
and the proof follows.
\end{proof}

\medskip
\noindent
\emph{Proof of Theorem~\ref{thm:main}}. 
Because of Lemma~\ref{lem:DF-recip}, part (b) follows 
by applying part (a) to $x^n h(1/x)$. We now prove 
part (a).

For the first statement, we only need to show that 
the polynomials
$a(x)$ and $b(x)$, defined by Equations~(\ref{eq:a}) 
and~(\ref{eq:b}), are real-rooted. This is clear for 
$b(x)$, since it is a nonnegative linear combination
of the elements of the interlacing sequence 
$\qQ_{\fF,n-1}$. By definition, we also have $a(x) = 
\sum_{i=0}^n \lambda_i(x) p_{\fF,n-1,i}(x)$ for some 
polynomials $\lambda_i(x)$ of degree at most one which 
have nonnegative coefficients and appear explicitly 
in~(\ref{eq:a}). The nonnegativity of the
$c_i$ easily implies that $(\lambda_n(x),
\lambda_{n-1}(x),\dots,\lambda_0(x))$ is an 
interlacing sequence. This observation, the fact 
(pointed out before Theorem~\ref{thm:Ath20+}) that 
$\pP_{\fF,n}$ is also interlacing and 
\cite[Lemma~7.8.3]{Bra15} imply that $a(x)$ is 
real-rooted as well.

For the second statement, let us write $h_\fF(x) = 
\dD_{\fF,n}(h(x))$. By \cite[Theorem~2.6]{BS20},
to prove that the real-rooted symmetric decomposition
of part (a) is interlacing, it suffices to show that
$h_\fF(x)$ is interlaced by $x^n h_\fF(1/x)$. For the 
latter, by \cite[Lemma~7.8.4]{Bra15}, it suffices to
show that $(\lambda x + \mu) x^n h_\fF(1/x) + h_\fF(x)$ 
is real-rooted for all positive reals $\lambda, \mu$. 
Since 
\[ x^n h_\fF(1/x) \ = \ \sum_{k=0}^n c_k x^n 
   p_{\fF,n,k}(1/x) \ = \ \sum_{k=0}^n c_k 
   p_{\fF,n,n-k}(x) , \]
we have 
\begin{eqnarray*} 
(\lambda x + \mu) x^n h_\fF(1/x) + h_\fF(x) & = & 
(\lambda x + \mu) \sum_{k=0}^n c_k p_{\fF,n,n-k}(x) 
\, + \, \sum_{k=0}^n c_k p_{\fF,n,k}(x) \\ & = & 
\sum_{k=0}^n \mu_k(x) p_{\fF,n,k}(x) , 
\end{eqnarray*}
where $\mu_k(x) = (c_k + c_{n-k} \mu) + c_{n-k} 
\lambda x$ for every $k \in \{0, 1,\dots,n\}$. Once 
again, it is routine to verify that the sequence 
$(\mu_n(x), \mu_{n-1}(x),\dots,\mu_0(x))$ is 
interlacing if $c_i c_{n-i-1} \le c_{i+1} c_{n-i}$ 
for every $i \in \{0, 1,\dots,n-1\}$. Another 
application of \cite[Lemma~7.8.3]{Bra15} then shows 
that $\sum_{k=0}^n \mu_k(x) p_{\fF,n,k}(x)$ is 
real-rooted and the proof follows.
\qed

\medskip
The following corollary produces large families 
of polynomials in face enumeration which admit 
nonnegative, real-rooted symmetric decompositions.
\begin{corollary} \label{cor:main}
Let $\fF$ be a feasible $f$-triangle which has the 
strong interlacing property with respect to $n$. 

\begin{itemize}
\itemsep=0pt
\item[{\rm (a)}] 
The polynomial $h_\fF(\Delta, x)$ has a nonnegative, 
real-rooted symmetric decomposition with respect to 
$n$ for every $(n-1)$-dimensional Cohen--Macaulay* 
simplicial complex $\Delta$. If, additionally, 
$\Delta$ satisfies the inequalities 
\begin{equation} \label{eq:tzanaki-ineq} 
\frac{h_0(\Delta)}{h_n(\Delta)} \, \le \, 
\frac{h_1(\Delta)}{h_{n-1}(\Delta)} \, \le \, 
\cdots \, \le \, 
\frac{h_{n-1}(\Delta)}{h_1(\Delta)} \, \le \, 
\frac{h_n(\Delta)}{h_0(\Delta)} , 
\end{equation}
then this decomposition is also interlacing.
 
\item[{\rm (b)}] 
The polynomial $h_\fF(\Delta, x)$ has a nonnegative, 
real-rooted symmetric decomposition with respect to 
$n-1$ for every triangulation $\Delta$ of the 
$(n-1)$-dimensional ball. If, additionally, $\Delta$ 
satisfies the inequalities 
\begin{equation} \label{eq:tzanakii-ineq} 
\frac{h_1(\Delta)}{h_{n-1}(\Delta)} \, \ge \, 
\frac{h_2(\Delta)}{h_{n-2}(\Delta)} \, \ge \, 
\cdots \, \ge \, 
\frac{h_{n-2}(\Delta)}{h_2(\Delta)} \, \ge \, 
\frac{h_{n-1}(\Delta)}{h_1(\Delta)} 
\end{equation}
(where terms involving an entry $h_i(\Delta) = 0$ 
may be ignored), then this decomposition is also 
interlacing.
\end{itemize}
\end{corollary}

\begin{proof}
This follows directly from Theorem~\ref{thm:main},
the fact that Cohen--Macaulay* simplicial complexes
satisfy (\ref{eq:ineqCM*}) and the fact (a consequence 
of \cite[Lemma~2.3]{Sta93}) that triangulations 
$\Delta$ of the $(n-1)$-dimensional ball satisfy the 
inequalities 
\[ h_0(\Delta) + h_1(\Delta)  + \cdots + h_i(\Delta) 
\ \ge \ h_{n-1}(\Delta) + h_{n-2}(\Delta) + \cdots + 
h_{n-i}(\Delta) \]
for all $1 \le i \le \lfloor n/2 \rfloor$.
\end{proof}

Let us record one situation in which the assumptions 
of Theorem~\ref{thm:main} on $h(x)$ are valid trivially.
\begin{corollary} \label{cor:small-degree}
Let $\fF$ be a feasible $f$-triangle which has the 
strong interlacing property with respect to $n$. 
For every $h(x) \in \RR[x]$ of degree at most $n/2$ 
with nonnegative coefficients, 
$\dD_{\fF,n}(h(x))$ has a nonnegative, real-rooted and 
interlacing symmetric decomposition with respect to 
$n-1$.

In particular, $h_\fF(\Delta, x)$ has such a 
decomposition for every $(n-1)$-dimensional simplicial 
complex $\Delta$ with nonnegative $h$-vector which 
satisfies $h_i(\Delta) = 0$ for $i \ge (n+1)/2$.
\end{corollary}

\begin{proof}
This follows from part (b) of Theorem~\ref{thm:main} 
since, under our assumptions, (\ref{eq:c-ineq2}) holds
trivially for $1 \le i \le \lfloor n/2 \rfloor$ and 
$c_{i+1} c_{n-i} = 0$ for every $i$.
\end{proof}

\section{Applications}
\label{sec:apps}

This section applies Theorem~\ref{thm:main} to the 
$r$-fold edgewise and $r$-colored barycentric 
subdivision and proves Theorems~\ref{thm:main-esdr} 
and~\ref{thm:main-sdr} and Corollary~\ref{cor:zono}.
We denote by $\fF$ and $\fF_{\esd_r}$ the 
$f$-triangles (of infinite size) defined by the 
barycentric and the $r$-fold edgewise subdivision,
respectively.

\subsection{The $r$-fold edgewise subdivision 
operator} \label{sec:esdr}

Recall from Section~\ref{sec:pre} that $\sS^r_k$ 
stands for the $k$th Veronese $r$-section operator. 
For the $r$-fold edgewise subdivision one has that
$\dD_{\fF_{\esd_r},n} = \uU^n_r: \RR_n[x] \to 
\RR_n[x]$, where
\begin{equation*} %%\label{eq:def-Unr}
\uU^n_r(h(x)) \ = \ \sS^r_0 \left( (1 + x + x^2 
+ \cdots + x^{r-1})^n h(x) \right) 
\end{equation*}
for $h(x) \in \RR_n[x]$; see, for instance, 
\cite[Section~4]{BW09} \cite[Section~4]{Ath14} 
\cite[Section~3]{Ath20+}. Equivalently,
\begin{eqnarray*}
\frac{\uU^n_r(h(x))}{(1-x)^n} & = & \frac{1} 
{(1-x)^n} \cdot \sS^r_0 
\left( (1 + x + x^2 + \cdots + x^{r-1})^n 
h(x) \right) \\ & & \\ & = & \sS^r_0 
\left( \frac{(1 + x + x^2 + \cdots + x^{r-1})^n} 
{(1-x^r)^n} h(x) \right) \ = \ \sS^r_0 
\left( \frac{h(x)}{(1-x)^n} \right),
\end{eqnarray*}
which shows that $\uU^n_r$ coincides with the 
operator which appears in the statement of 
Theorem~\ref{thm:main-esdr} under the same name.

\medskip
\noindent
\emph{Proof of Theorem~\ref{thm:main-esdr}}. 
This follows directly from Theorem~\ref{thm:main} 
and the fact \cite[Example~7.2]{Ath20+} that the 
$f$-triangle of the $r$-fold edgewise subdivision 
has the strong interlacing property with respect 
to $n$ for every $r \ge n$.
\qed

\medskip
The conclusion of part (b) of 
Theorem~\ref{thm:main-esdr} was proven in \cite{Jo20+}
under stronger assumptions (see \cite[Theorem~1.1]{Jo20+}) 
which, for example, do not cover part (b) of the 
following corollary.
\begin{corollary} \label{cor:esdr}
\begin{itemize}
\itemsep=0pt
\item[{\rm (a)}] 
The polynomial $h(\esd_r(\Delta), x)$ has a nonnegative, 
real-rooted symmetric decomposition with respect to 
$n$ for every $r \ge n$ and every $(n-1)$-dimensional 
Cohen--Macaulay* simplicial complex $\Delta$. 
 
\item[{\rm (b)}] 
The polynomial $h(\esd_r(\Delta), x)$ has a nonnegative, 
real-rooted symmetric decomposition with respect to 
$n-1$ for every $r \ge n$ and every triangulation 
$\Delta$ of the ball of dimension $n-1$. 
\end{itemize}
\end{corollary}

\begin{proof}
Apply Theorem~\ref{thm:main-esdr} to $h(\Delta,x)$ or, 
alternatively, Corollary~\ref{cor:main} to the 
$r$-fold edgewise subdivision.
\end{proof}

The following statement improves 
\cite[Proposition~5.2]{Jo20+}.
\begin{corollary} \label{cor:esdr-small}
The polynomial $\uU^n_r(h(x))$ has a nonnegative, 
real-rooted and interlacing symmetric decomposition 
with respect to $n-1$ for every 
polynomial $h(x) \in \RR[x]$ of degree at most $n/2$ 
with nonnegative coefficients and every $r \ge n$.
\end{corollary}

\begin{proof}
This follows from Corollary~\ref{cor:small-degree} 
and the fact that the $f$-triangle of the $r$-fold 
edgewise subdivision has the strong interlacing 
property with respect to $n$ for every $r \ge n$. 
\end{proof}

\begin{figure}
\begin{center}
\begin{tikzpicture}[scale=3]
\label{fg:sd3}
			
			\coordinate (a) at (1,0);
			\coordinate (b) at (0,1.73);
			\coordinate (c) at (-1,0);
			\coordinate (o) at (0,0.58);
			
			\coordinate (ab) at ($(a)!0.5!(b)$);
			\coordinate (bc) at ($(b)!0.5!(c)$);
			\coordinate (ca) at ($(c)!0.5!(a)$);
			\coordinate (b1) at ($(b)!0.33!(o)$);
			\coordinate (b2) at ($(b)!0.66!(o)$);
			\coordinate (b3) at ($(o)!0.33!(ca)$);
			\coordinate (b4) at ($(o)!0.66!(ca)$);
			
			\coordinate (c1) at ($(c)!0.33!(o)$);
			\coordinate (c2) at ($(c)!0.66!(o)$);
			\coordinate (c3) at ($(o)!0.33!(ab)$);
			\coordinate (c4) at ($(o)!0.66!(ab)$);	
			
			\coordinate (a1) at ($(a)!0.33!(o)$);
			\coordinate (a2) at ($(a)!0.66!(o)$);
			\coordinate (a3) at ($(o)!0.33!(bc)$);
			\coordinate (a4) at ($(o)!0.66!(bc)$);	
			
			\coordinate (ab1) at ($(a)!0.16!(b)$);	
			\coordinate (ab2) at ($(a)!0.33!(b)$);	
			\coordinate (ab4) at ($(a)!0.66!(b)$);	
			\coordinate (ab5) at ($(a)!0.83!(b)$);	
			
			\coordinate (bc1) at ($(b)!0.16!(c)$);	
			\coordinate (bc2) at ($(b)!0.33!(c)$);	
			\coordinate (bc4) at ($(b)!0.66!(c)$);	
			\coordinate (bc5) at ($(b)!0.83!(c)$);	
			
			\coordinate (ca1) at ($(c)!0.16!(a)$);	
			\coordinate (ca2) at ($(c)!0.33!(a)$);	
			\coordinate (ca4) at ($(c)!0.66!(a)$);	
			\coordinate (ca5) at ($(c)!0.83!(a)$);	
			
			\draw (a)--(b)--(c)--(a);
			\draw (c)--(ab);
			\draw (a)--(bc);
			\draw (b)--(ca);
			
			\draw[dashed] (bc1)--(b1)--(ab5);
			\draw[dashed] (bc2)--(b2)--(ab4);
			\draw[dashed] (bc4)--(c2)--(ca2);
			\draw[dashed] (bc5)--(c1)--(ca1);
			
			\draw[dashed] (ab2)--(a2)--(ca4);
			\draw[dashed] (ab1)--(a1)--(ca5);
			
			\draw[dashed] (bc1)--(a3)--(bc5);
			\draw[dashed] (bc2)--(a4)--(bc4);
			\draw[dashed] (ab5)--(c3)--(ab1);
			\draw[dashed] (ab4)--(c4)--(ab2);
			
			\draw[dashed] (ca1)--(b3)--(ca5);
			\draw[dashed] (ca2)--(b4)--(ca4);
			
			\draw[dashed] (c1)--(b1);
			\draw[dashed] (c2)--(b2);
			\draw[dashed] (c1)--(a1)--(b1);
			\draw[dashed] (c2)--(a2)--(b2);
			
			\fill (a)circle (0.75pt);
			\fill (b)circle (0.75pt);
			\fill (c)circle (0.75pt);
			\fill (ab)circle (0.75pt);
			\fill (bc)circle (0.75pt);
			\fill (ca)circle (0.75pt);
			\fill (o)circle (0.75pt);
			\fill (b1)circle (0.5pt);
			\fill (b2)circle (0.5pt);
			\fill (b3)circle (0.5pt);
			\fill (b4)circle (0.5pt);
			\fill (c1)circle (0.5pt);
			\fill (c2)circle (0.5pt);
			\fill (c3)circle (0.5pt);
			\fill (c4)circle (0.5pt);
			\fill (a1)circle (0.5pt);
			\fill (a2)circle (0.5pt);
			\fill (a3)circle (0.5pt);
			\fill (a4)circle (0.5pt);
			\fill (ab1)circle (0.5pt);
			\fill (ab2)circle (0.5pt);
			\fill (ab4)circle (0.5pt);
			\fill (ab5)circle (0.5pt);
			
			\fill (bc1)circle (0.5pt);
			\fill (bc2)circle (0.5pt);
			\fill (bc4)circle (0.5pt);
			\fill (bc5)circle (0.5pt);
			
			\fill (ca1)circle (0.5pt);
			\fill (ca2)circle (0.5pt);
			\fill (ca4)circle (0.5pt);
			\fill (ca5)circle (0.5pt);
						
\end{tikzpicture}
\caption{The 3-colored barycentric subdivision of the 
2-simplex}
\end{center}
\end{figure}
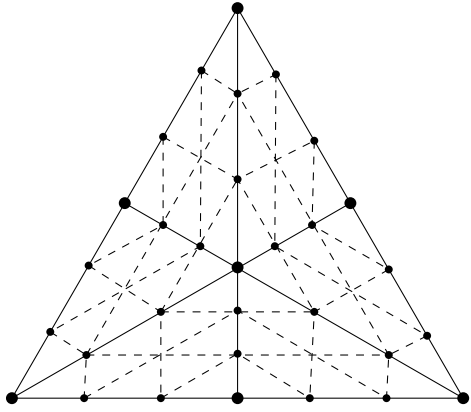

\subsection{The $r$-colored barycentric subdivision 
operator}
\label{sec:r-colored}

Consider the composition of linear operators 
$\dD_{n,r} := \dD_{\fF_{\esd_r},n} \circ \dD_{\fF,n} 
= \uU^n_r \circ \dD_n: \RR_n[x] \to \RR_n[x]$. 
Thus,
\begin{equation} \label{eq:h-sdr-alt}
\dD_{n,r}(h(x)) \ = \ \uU^n_r(\dD_n(h(x))) \ = \ 
\sS^r_0 \left( (1 + x + x^2 + \cdots + x^{r-1})^n 
\dD_n(h(x)) \right)
\end{equation}
for every $h(x) \in \RR_n[x]$. We confirm in 
the proof of Theorem~\ref{thm:main-sdr}, given in 
the sequel, that $\dD_{n,r}$ coincides with the 
operator which appears in the statement of this
theorem, under the same name. 

Clearly, we have $\dD_{n,r} = \dD_{\fF_r,n}$, 
where $\fF_r$ is the $f$-triangle of the uniform 
triangulation obtained by first taking the 
barycentric subdivision of a simplicial complex 
$\Delta$ and then the $r$-fold edgewise subdivision 
of that. This triangulation (see Figure~1 for an 
example), termed as the \emph{$r$-colored 
barycentric subdivision} 
in~\cite{Ath20+}, was introduced in \cite{Ath14} 
in order to partially interpret geometrically the 
derangement polynomial for the colored permutation 
group $\ZZ_r \wr \fS_n$; it was further studied 
enumeratively in~\cite{Ath20b}. As mentioned 
in~\cite[Section~3]{Ath20+}, $\fF_2$ coincides with
the $f$-triangle defined by the interval 
triangulation \cite[Section~3.3]{MW17}, so all our 
results here apply to that as well.

The main result of this section answers in the 
affirmative a question from~\cite[Section~7]{Ath20+}.
\begin{theorem} \label{thm:sdr-strong} 
The $f$-triangle $\fF_r$ has the strong interlacing 
property.
\end{theorem}

We postpone the proof until the end of the section.

\medskip
\noindent
\emph{Proof of Theorem~\ref{thm:main-sdr}}. This 
follows from Theorems~\ref{thm:Ath20+},~\ref{thm:main} 
and~\ref{thm:sdr-strong}, 
provided $\dD_{n,r}$ coincides with the operator
in the statement of Theorem~\ref{thm:main-sdr}.
Indeed, considering $h(x) \in \RR_n[x]$ and setting 
$f(x) = \sum_{i=0}^n c_i x^i(1+x)^{n-i}$, from 
Equation~(\ref{eq:h-sdr-alt}) we get 

\begin{eqnarray*}
\frac{\dD_{n,r}(h(x))}{(1-x)^n} & = & \frac{1} 
{(1-x)^n} \cdot \sS^r_0 
\left( (1 + x + x^2 + \cdots + x^{r-1})^n 
\dD_n(h(x)) \right) \\ & & \\ & = & \sS^r_0 
\left( \frac{(1 + x + x^2 + \cdots + x^{r-1})^n} 
{(1-x^r)^n} \dD_n(h(x)) \right) \ = \ 
\sS^r_0 \left( \frac{\dD_n(h(x))}{(1-x)^n} \right) 
\\ & & \\ & = & \sS^r_0 \left( (1-x) 
\sum_{m \ge 0} f(m) x^m  \right) \ = \ 
f(0) + \sum_{m \ge 1} \left( f(rm) - f(rm-1) 
\right) x^m
\end{eqnarray*}

\medskip
\noindent
and the proof follows. Note that we have used 
Equation~(\ref{eq:h-sd}) for the next to last step.
\qed

\medskip
Part (a) of Theorem~\ref{thm:main-sdr} was deduced 
in \cite[Proposition~7.5]{Ath20+} from 
Theorem~\ref{thm:B-BS} and the fact that $\uU^n_r$ 
preserves real-rootedness for polynomials with 
nonnegative coefficients. The present proof shows 
additionally that $\dD_{n,r}(h(x))$ is interlaced 
by $h_{\fF_r}
(\sigma_n, x)$ and interlaces $x^n h_{\fF_r}(\sigma_n,
1/x)$; see \cite[Proposition~5.1]{Ath20b}
for combinatorial interpretations of $h_{\fF_r}
(\sigma_n,x)$.

Part (b) of the following corollary can be deduced 
from \cite[Theorem~1.1]{Jo20+} for $r \ge n-1$ but,
to the best of our knowledge, not for other values 
of $r$.
\begin{corollary} \label{cor:sdr}
Let $\sd_r(\Delta)$ denote the $r$-colored 
barycentric subdivision of $\Delta$.

\begin{itemize}
\itemsep=0pt
\item[{\rm (a)}] 
The polynomial $h(\sd_r(\Delta), x)$ has a nonnegative, 
real-rooted symmetric decomposition with respect to 
$n$ for every $(n-1)$-dimensional Cohen--Macaulay* 
simplicial complex $\Delta$. 
 
\item[{\rm (b)}] 
The polynomial $h(\sd_r(\Delta), x)$ has a nonnegative, 
real-rooted symmetric decomposition with respect to 
$n-1$ for every triangulation $\Delta$ of the 
$(n-1)$-dimensional ball. 
\end{itemize}
\end{corollary}

\begin{proof}
Apply Theorem~\ref{thm:main-sdr} to $h(\Delta,x)$ or, 
alternatively, Corollary~\ref{cor:main} to the 
$r$-colored barycentric subdivision.
\end{proof}

\noindent
\emph{Proof of Corollary~\ref{cor:zono}}. It was 
observed in \cite[Section~4]{BS20} that $\iota(\zZ; 
x) = \sum_{i=0}^n c_i x^i(1+x)^{n-i}$ for some 
nonnegative integers $c_i$ (these are the coefficients
of the $h^*$-polynomial of the Lawrence polytope 
associated to $\zZ$) which satisfy inequalities
(\ref{eq:c-ineq1}) for all $0 \le i \le \lfloor n/2 
\rfloor$. Hence, the result follows by applying 
Theorem~\ref{thm:main-sdr} to $f(x) = \iota(\zZ; 
x)$.
\qed

\medskip
We now turn our attention to the proof of 
Theorem~\ref{thm:sdr-strong}. Recall that $\fF = 
\fF_1$ is the $f$-triangle defined by barycentric 
subdivision. We write $p_{n,k}(x)$ in place of 
$p_{\fF,n,k}(x)$ and recall that 
\begin{equation}
\label{eq:pnkrec-long}
p_{n,k}(x) \ = \ x \sum_{i=0}^{k-1} 
p_{n-1,i}(x) \, + \, \sum_{i=k}^{n-1} p_{n-1,i}(x) 
\end{equation}
for all $n \ge 1$ and $k \in \{0, 1,\dots,n\}$ since, 
as mentioned in Section~\ref{sec:uniform}, 
$p_{\fF,n-1,n}(x) = 0$ in this case. A combinatorial 
interpretation of the polynomials $p_{\fF_r,n,k}(x)$ 
for $k \in \{0, 1,\dots,n\}$ was given in
\cite[Proposition~4.7]{Ath20+}. We will now express 
these polynomials and $\theta_{\fF_r}(\sigma_n, x)$ 
in terms of the $p_{n,k}(x)$. For this reason, we 
introduce the polynomials
\begin{equation*} %%\label{eq:def-pnkri}
p^{\langle r,j \rangle}_{n,k}(x) \ = \ \sS^r_j 
\left( (1 + x + x^2 + \cdots + x^{r-1})^n 
p_{n,k}(x) \right) 
\end{equation*}
for $n \in \NN$, $k \in \{0, 1,\dots,n\}$ and $j \in 
\{0, 1,\dots,r-1\}$. For $r=2$ they have been considered 
before in~\cite{AN20} (see, for instance, Corollary~4.7 
there).

\bigskip
{\scriptsize
\begin{table}[hptb]
\begin{center}
\begin{tabular}{| l || l | l | l | l | l ||} \hline
& $k=0$ & $k=1$ & $k=2$    & $k=3$ \\ \hline \hline
$j=0$   & $1+34x+19x^2$ & $30x+24x^2$ & $24x+30x^2$ & 
          $19x+34x^2+x^3$ \\ \hline
$j=1$   & $7+40x+7x^2$ & $4+40x+10x^2$ & $2+38x+14x^2$ 
        & $1+34x+19x^2$ \\ \hline
$j=2$   & $19+34x+x^2$ & $14+38x+2x^2$ & 
          $10+40x+4x^2$ & $7+40x+7x^2$ \\ \hline
\end{tabular}
\caption{The polynomials 
$p^{\langle r,j \rangle}_{n,k}(x)$ for $n=r=3$.}
\label{tab:pnkrj}
\end{center}
\end{table}}

\begin{proposition} \label{prop:sdr-ptheta} 
For the $f$-triangle $\fF_r$ we have
\begin{equation} \label{eq:Fr-pnk}
p_{\fF_r,n,k}(x) \ = \ 
         p^{\langle r,0 \rangle}_{n,k}(x)
\end{equation}
for all $n \in \NN$ and $k \in \{0, 1,\dots,n\}$ 
and
\begin{equation} \label{eq:Fr-theta}
\theta_{\fF_r}(\sigma_n, x) \ = \ x \, \sum_{j=1}^{r-1} 
\sum_{k=0}^{n-1} p^{\langle r,j \rangle}_{n-1,k}(x)
\end{equation}
for every $n \ge 1$.
\end{proposition}

\begin{proposition} \label{prop:pnkrj-rec} 
The polynomials $p^{\langle r,j \rangle}_{n,k}(x)$ 
satisfy the recurrence
\[ p^{\langle r,j \rangle}_{n,k}(x) \ = \ 
x \sum_{\ell=j+1}^{r-1} \sum_{i=0}^{n-1} 
p^{\langle r, \ell \rangle}_{n-1,i}(x) \, + \, 
x \sum_{i=0}^{k-1} 
p^{\langle r,j \rangle}_{n-1,i}(x) \, + \, 
\sum_{i=k}^{n-1} 
p^{\langle r,j \rangle}_{n-1,i}(x) \, + \, 
\sum_{\ell=0}^{j-1} \sum_{i=0}^{n-1} 
p^{\langle r, \ell \rangle}_{n-1,i}(x) \]
for every $n \ge 1$ and all $k \in \{0, 1,\dots,n\}$ 
and $j \in \{0, 1,\dots,r-1\}$.
\end{proposition}

\begin{proof}
From the definition of $p^{\langle r,j \rangle}_{n,k}
(x)$ we get 

\begin{eqnarray*}
p^{\langle r,j \rangle}_{n,k}(x) & = &
\sS^r_j \left( (1 + x + x^2 + \cdots + x^{r-1})^n 
p_{n,k}(x) \right) \\ & & \\ & = & 
\sS^r_j \left( \left( \, \sum_{\ell=0}^{r-1} x^\ell 
\right) (1 + x + x^2 + \cdots + x^{r-1})^{n-1} 
p_{n,k}(x) \right) \\ & = & 
\sum_{\ell=0}^{r-1} \sS^r_j \left( x^\ell 
(1 + x + x^2 + \cdots + x^{r-1})^{n-1} p_{n,k}(x) 
\right) .
\end{eqnarray*}

\medskip
\noindent
Replacing $p_{n,k}(x)$ by the right-hand side of 
(\ref{eq:pnkrec-long}) and changing the order of 
summation, we get 

\begin{eqnarray*}
p^{\langle r,j \rangle}_{n,k}(x) & = &
\sum_{i=0}^{k-1} \sum_{\ell=0}^{r-1} \sS^r_j 
\left( x^{\ell+1} (1 + x + x^2 + \cdots + x^{r-1})^{n-1} 
p_{n-1,i}(x) \right) + \\ &  & 
\sum_{i=k}^{n-1} \sum_{\ell=0}^{r-1} \sS^r_j \left( 
x^\ell (1 + x + x^2 + \cdots + x^{r-1})^{n-1} 
p_{n-1,i}(x) \right) 
\end{eqnarray*}

\medskip
\noindent
and applying (\ref{eq:Sx}) yields the desired expression 
for $p^{\langle r,j \rangle}_{n,k}(x)$; the details 
are omitted.
\end{proof}

\medskip
\noindent
\emph{Proof of Proposition~\ref{prop:sdr-ptheta}}. 
For every $h(x) = c_0 + c_1 x + \cdots + c_n x^n 
\in \RR_n[x]$,

\begin{eqnarray*}
\dD_{\fF_r,n} (h(x)) & = & 
\uU^n_r (\dD_n (h(x)) \ = \ \uU^n_r 
\left( \, \sum_{k=0}^n c_k p_{n,k}(x) \right) \ = \
\sum_{k=0}^n c_k \, \uU^n_r \left( p_{n,k}(x) \right) 
\\ & = & \sum_{k=0}^n c_k \, \sS^r_0 \left( 
(1 + x + x^2 + \cdots + x^{r-1})^n p_{n,k}(x) \right) \ 
= \ \sum_{k=0}^n c_k p^{\langle r,0 \rangle}_{n,k}(x) .
\end{eqnarray*}

\medskip
\noindent
This proves~(\ref{eq:Fr-pnk}). Since $h_{\fF_r} 
(\sigma_n, x) = p_{\fF_r,n,0}(x) = 
p^{\langle r,0 \rangle}_{n,0}(x)$, from the recurrence 
of Proposition~\ref{prop:pnkrj-rec} we get
\[ h_{\fF_r} (\sigma_n, x) \ = \  x \sum_{j=1}^{r-1} 
\sum_{k=0}^{n-1} p^{\langle r,r-j \rangle}_{n-1,k}(x) 
\, + \, \sum_{k=0}^{n-1} 
p^{\langle r,0 \rangle}_{n-1,k}(x) . \]
Since $\theta_{\fF_r}(\sigma_n, x) = h_{\fF_r}
(\sigma_n, x) - h_{\fF_r}(\partial \sigma_n, x)$ and 
\[ h_{\fF_r}(\partial \sigma_n, x) \ = \ \sum_{k=0}^{n-1} 
p_{\fF_r,n-1,k}(x) \ = \ \sum_{k=0}^{n-1} 
p^{\langle r,0 \rangle}_{n-1,k}(x)  \]
by (\ref{eq:pnk-sum}) and part (a), the proof 
of~(\ref{eq:Fr-theta}) follows. 
\qed

\medskip
\noindent
\emph{Proof of Theorem~\ref{thm:sdr-strong}}. For 
$j \in \{0, 1,\dots,r-1\}$ we consider the sequence 
\[ \pP^{\langle r,j \rangle}_n \ := \ 
(p^{\langle r,j \rangle}_{n,0}(x))_{0 \le k \le n} \ 
= \ (p^{\langle r,j \rangle}_{n,0}(x), 
p^{\langle r,j \rangle}_{n,1}(x),\dots,p^{\langle r,j 
\rangle}_{n,n}(x))  \]
and let 
\[ \pP_{n,r} \ := \  
   (\pP^{\langle r,r-j \rangle}_n)_{1 \le j \le r} 
	 \ = \ (\pP^{\langle r,r-1 
   \rangle}_n,\dots,\pP^{\langle r,1 \rangle}_n,
   \pP^{\langle r,0 \rangle}_n) \]
be their concatenation, in the specified order; see
see Table~\ref{tab:pnkrj} for an example. 

We claim that $\pP_{n,r}$ is interlacing
for every $n \in \NN$. This is clear for $n=0$, since 
$\pP_{0,r} = (0,\dots,0, 1)$, so we assume that $n \ge 
1$. Proposition~\ref{prop:pnkrj-rec} implies that 
$p^{\langle r,j \rangle}_{n,n}(x) = p^{\langle r,j-1 
\rangle}_{n,0}(x)$ for every $j \in \{1, 2,\dots,r-1\}$, 
so $\pP_{n,r}$ has 
$r-1$ pairs of equal consecutive elements. The same 
proposition shows that, when one of these elements is 
removed from each of these pairs, the resulting 
sequence is obtained from $\pP_{n-1,r}$ by the recipe 
of~(\ref{eq:recipe}). As a result, and since 
doubling some elements of an interlacing sequence
clearly preserves the interlacing property, the 
interlacing of $\pP_{n-1,r}$ implies that of 
$\pP_{n,r}$ and our claim follows by induction 
on $n$.

We may now prove the theorem. Clearly, the 
polynomials $p^{\langle r,j \rangle}_{n,k}(x)$ 
have nonnegative coefficients. As can be inferred 
from their definition or 
Proposition~\ref{prop:pnkrj-rec}, they have degree 
$n-1$ except for $p^{\langle r,0 \rangle}_{n,n}
(x)$, which has degree $n$. Thus,
Proposition~\ref{prop:sdr-ptheta} shows that 
$\theta_{\fF_r}(\sigma_n, x)$ has nonnegative 
coefficients and degree $n-1$. Given that 
$\pP_{n-1,r}$ is an interlacing sequence, it 
also shows that $\theta_{\fF_r}(\sigma_n, x)/x$ 
is a sum of polynomials each of which interlaces 
$p^{\langle r,0 \rangle}_{n-1,0}(x) = h_{\fF_r}
(\sigma_{n-1}, x)$. This implies that 
$\theta_{\fF_r}(\sigma_n, x)$ is real-rooted and 
interlaced by $h_{\fF_r}(\sigma_{n-1}, x)$ and 
the proof follows.
\qed

\section{Skeleta of simplicial complexes}
\label{sec:skeleta}

This section proves and generalizes  
Theorem~\ref{thm:skel-intro} in the setting of uniform 
triangulations as follows.
\begin{theorem} \label{thm:skel} 
Let $\Gamma$ be an $n$-dimensional simplicial 
complex with nonnegative $h$-vector and let $\Delta$ 
be the $(n-1)$-dimensional skeleton of $\Gamma$. 

\begin{itemize}
\itemsep=0pt
\item[{\rm (a)}]
The polynomial $h_\fF(\Delta, x)$ has a nonnegative, 
real-rooted and interlacing symmetric decomposition 
with respect to $n$ for every feasible $f$-triangle 
$\fF$ which has the strong interlacing property with 
respect to $n$. 

\item[{\rm (b)}] 
The polynomial $h_\fF(\Delta, x)$ interlaces $h_\fF
(\Gamma, x)$ for every feasible $f$-triangle $\fF$ 
which has the strong interlacing property with 
respect to $n+1$.
\end{itemize}
\end{theorem}

\begin{proof}
As a direct consequence of the defining equation 
(\ref{eq:hdef}) of the $h$-polynomial, the entries 
of the $h$-vector of $\Delta$ can be expressed in 
terms of those of the $h$-vector of $\Gamma$ as 
\[ h_k(\Delta) \ = \ h_0(\Gamma) + h_1(\Gamma) + 
                     \cdots + h_k(\Gamma) \]
for $0 \le k \le n$. In particular, $h_0(\Delta) 
\le h_1(\Delta) \le \cdots \le h_n(\Delta)$ and 
this makes it obvious that the $h_i(\Delta)$ 
satisfy the inequalities~(\ref{eq:c-ineq1}) and 
that $h_i(\Delta) h_{n-i-1}(\Delta) \le h_{i+1}
(\Delta)h_{n-i}(\Delta)$ for $0 \le i \le n-1$. 
Thus, part (a) follows from part (a) of 
Theorem~\ref{thm:main}.

Since $\fF$ has the strong interlacing property 
with respect to $n+1$, $h_\fF(\Delta, x)$ and 
$h_\fF(\Gamma, x)$ have nonnegative coefficients 
and only real roots by Theorem~\ref{thm:Ath20+}. 
Moreover, by Equation~(\ref{eq:hF}) and our previous 
remarks,
\begin{eqnarray*}
h_\fF(\Delta, x) & = & \sum_{k=0}^n \left( 
h_0(\Gamma) + h_1(\Gamma) + \cdots + h_k(\Gamma) 
\right) p_{\fF,n,k}(x) \\ & & \\
h_\fF(\Gamma, x) & = & \sum_{k=0}^{n+1} h_k(\Gamma) 
\, p_{\fF,n+1,k}(x) \, . 
\end{eqnarray*}
Expressing the $p_{\fF,n+1,k}(x)$ in terms of the 
$p_{\fF,n,k}(x)$ by (\ref{eq:pnk-rec-long}), we 
compute that for any positive reals $\lambda, \mu$,
\[ (\lambda x + \mu) h_\fF(\Delta, x) + h_\fF
   (\Gamma, x) \ = \ \sum_{k=0}^{n+1} \nu_k(x) 
   p_{\fF,n,k}(x) \, , \]
where 
\[ \nu_k(x) \ = \ \begin{cases}
  {\displaystyle (\mu+1) \sum_{i=0}^k h_i(\Gamma) + \left( 
	\lambda \sum_{i=0}^k h_i(\Gamma) + \sum_{i=k+1}^{n+1} 
	h_i(\Gamma) \right) x}, & \text{if $0 \le k \le n$} \\
  {\displaystyle \sum_{i=0}^{n+1} h_i(\Gamma)}, & 
	\text{if $k = n+1$.} \end{cases} \]

\medskip
\noindent
As in the proof of Theorem~\ref{thm:main}, it is 
routine to show that $(\nu_{n+1}(x),\dots,\nu_1(x), 
\nu_0(x))$ is an interlacing sequence. Since 
$\pP_{\fF,n+1}$ is also interlacing by the proof of 
\cite[Theorem~6.1]{Ath20+}, the result of part (b) 
follows by applying Lemmas~7.8.3 
and~7.8.4 of~\cite{Bra15}.
\end{proof}

%\medskip
\noindent
\emph{Proof of Theorem~\ref{thm:skel-intro}}. 
Apply Theorem~\ref{thm:skel} to the $f$-triangles 
of the $r$-fold edgewise and $r$-colored 
barycentric subdivisions.
\qed 

\section{Concluding remarks and open problems}
\label{sec:rem}

Given the crucial role played by the strong 
interlacing property in Theorems~\ref{thm:Ath20+} 
and~\ref{thm:main}, the following question arises
naturally. 
\begin{question} \label{que:ath20}
Which uniform triangulations have the strong 
interlacing property?
\end{question}

The inequalities (\ref{eq:tzanaki-ineq}) imply
that $h_i(\Delta) \le h_{n-i}(\Delta)$ for all
$0 \le i \le \lfloor n/2 \rfloor$. The validity 
of the latter inequalities for doubly 
Cohen--Macaulay complexes follows from 
\cite[Corollary~6.2]{APP21} (and was earlier 
shown for the more restrictive class of simplicial 
complexes with a convex ear decomposition in
\cite[Corollary~3.10]{Sw06}). Similar remarks 
apply to inequalities (\ref{eq:tzanakii-ineq}). 
We thus ask the following questions. 
\begin{question} \label{que:tzanaki-ineq}
Which Cohen--Macaulay* simplicial complexes 
satisfy (\ref{eq:tzanaki-ineq})? Do these 
inequalities hold for all doubly Cohen--Macaulay 
simplicial complexes?
\end{question}
\begin{question} \label{que:tzanakii-ineq}
Which triangulations of the ball satisfy
(\ref{eq:tzanakii-ineq})? 
\end{question}

We expect that Question~\ref{que:tzanaki-ineq} 
has an affirmative answer at least for interesting 
classes of doubly Cohen--Macaulay complexes. The 
inequalities $h_i(\Delta) \le h_{n-i}(\Delta)$ we
mentioned earlier, the fact that doubly Cohen--Macaulay 
simplicial complexes are level \cite[p.~94]{StaCCA}
and \cite[Proposition~III.3.3~(a)]{StaCCA} imply 
an affirmative answer for doubly Cohen--Macaulay 
complexes of dimension at most 3. 
Question~\ref{que:tzanakii-ineq} has an affirmative
answer in three dimensions as well, since $h_1(\Delta) 
\ge h_3(\Delta)$ for every triangulation $\Delta$ of 
the 3-dimensional ball (see, for instance, 
\cite[Section~3]{Ko11}) but not for every 
triangulated ball $\Delta$ in higher dimensions. 
Indeed, according to \cite[Theorem~14]{Ko11},
$(1, a, b, 1, 1)$ is the $h$-vector of a 
triangulation of the 4-dimensional ball for all 
positive integers $a, b$ with $b \le 1 + a(a-1)/2$.

Some of the problems about simplicial complexes we 
have studied make sense for polyhedral (or even more 
general cell) complexes. We record two of them here, 
one of which has already been mentioned in the 
introduction. Let $\sd(\lL)$ denote the 
barycentric subdivision of a (finite) polyhedral 
complex $\lL$.
\begin{question} 
\begin{itemize}
\itemsep=0pt
\item[{\rm (a)}]
Does $h(\sd(\lL),x)$ have a nonnegative, real-rooted 
symmetric decomposition with respect to $n$ for every 
$(n-1)$-dimensional Cohen--Macaulay* polyhedral 
complex $\lL$? 

\item[{\rm (b)}] 
Does $h(\sd(\lL),x)$ have a nonnegative, real-rooted 
symmetric decomposition with respect to $n-1$ for every 
$(n-1)$-dimensional polyhedral ball $\lL$? 
\end{itemize}

If so, are these decompositions interlacing? 
\end{question}

Finally, we noticed in Section~\ref{sec:skeleta} 
that the $(n-1)$-skeleton of any $n$-dimensional 
Cohen--Macaulay simplicial complex has an 
increasing $h$-vector.
\begin{question}
Which Cohen--Macaulay simplicial complexes have 
increasing $h$-vector?
\end{question}

%% \bigskip
%% \noindent \textbf{Acknowledgments}. 

\end{document}